\documentclass[11pt]{article}
\usepackage[margin=1in]{geometry}
\usepackage{amsmath,amssymb,amsthm}
\usepackage{hyperref}
\usepackage[nameinlink,capitalize,noabbrev]{cleveref}

\usepackage{titlesec}
\titleformat{\section}
{\normalfont\fontsize{12}{15}\bfseries}{\thesection}{1em.}{}

\usepackage[labelfont=bf]{caption}
\counterwithin{figure}{section}
\counterwithin{table}{section}

\newtheorem{theorem}{Theorem}[section]
\newtheorem{proposition}[theorem]{Proposition}

\newtheorem{corollary}[theorem]{Corollary}
\newtheorem{lemma}[theorem]{Lemma}
\newtheorem{definition}[theorem]{Definition}
\newtheorem{remark}[theorem]{Remark}
\newtheorem{example}[theorem]{Example}

\title{Conditional Park--Pham Bounds under Positive Correlation}
\author{
Bryce Alan Christopherson\thanks{Department of Mathematics and Statistics, University of North Dakota, Grand Forks, ND, USA. Email: bryce.christopherson@UND.edu}
\and
Darian Colgrove\thanks{Department of Mathematics and Statistics, University of North Dakota, Grand Forks, ND, USA.}
}
\date{}

\begin{document}

\maketitle 

\normalsize

 \begin{abstract}
 \noindent
        We record a conditional form of the $\epsilon$-dependent Park--Pham theorem.  If a monotone property $\mathcal{F}\subseteq 2^X$ is positively correlated with a conditioning event $B\subseteq 2^X$ under the product measure $\mu_p$, then the usual Park--Pham lower bound for $\mu_p(\mathcal{F})$ transfers to the conditional probability $\mathbb P(X_p\in\mathcal{F}\mid X_p\in B)$.  This gives, in particular, conditional Park--Pham bounds for increasing conditioning events by Harris's inequality, and for nonmonotone conditioning events that are independent of the target property.  We also formulate the transfer principle for finite posets embedded in Boolean lattices and illustrate it with pattern-containment upper sets in permutation classes. \\[2mm]
 \end{abstract}

\baselineskip=0.20in

\section{Introduction and preliminaries}
    
    Threshold phenomena for monotone properties form a central theme in probabilistic combinatorics \cite{park-threshold-phenomena}.  Let $X$ be a finite set, and write $2^X$ for the Boolean lattice ordered by inclusion.  A family $\mathcal{F}\subseteq 2^X$ is called an \textit{upper set}, or increasing event, if $S\in\mathcal{F}$ and $S\subseteq T\subseteq X$ implies $T\in\mathcal{F}$.  We say that $\mathcal{F}$ is nontrivial if $\mathcal{F}\neq\emptyset,2^X$.
    
    For $p\in[0,1]$, let $\mu_p$ be the product measure on $2^X$ defined by $\mu_p(S)=p^{|S|}(1-p)^{|X|-|S|}$ for $S\subseteq X$.  For a family $\mathcal A\subseteq 2^X$, write $\mu_p(\mathcal A)=\sum_{S\in\mathcal A}\mu_p(S)$.  Equivalently, if $X_p$ is a random subset of $X$ with distribution $\mu_p$, then $\mu_p(\mathcal A)=\mathbb P(X_p\in\mathcal A)$.  For a nontrivial upper set $\mathcal{F}$, the \textit{critical probability} $p_c(\mathcal{F})$ is the unique value of $p$ for which $\mu_p(\mathcal{F})=\frac{1}{2}$.
    
    The Park--Pham theorem \cite{parkpham,shortparkphamproof}, formerly the Kahn-Kalai conjecture \cite{kahn-kalai}, bounds $p_c(\mathcal{F})$ in terms of the expectation threshold $q(\mathcal{F})$.  If $\mathcal G\subseteq 2^X$, write $\langle \mathcal G\rangle=\{T\subseteq X:S\subseteq T\text{ for some }S\in\mathcal G\}$.  A family $\mathcal G$ is a cover of $\mathcal{F}$ if $\mathcal{F}\subseteq \langle\mathcal G\rangle$.  The upper set $\mathcal{F}$ is called $p$-small if there exists a cover $\mathcal G$ of $\mathcal{F}$ such that $\sum_{S\in\mathcal G}p^{|S|}\leq \frac{1}{2}$.  The \textit{expectation threshold} $q(\mathcal{F})$ is the supremum of the values $p$ for which $\mathcal{F}$ is $p$-small.  Let $\mathcal{F}_0$ denote the set of minimal elements of $\mathcal{F}$, and write $\ell(\mathcal{F})=\max\{|S|:S\in\mathcal{F}_0\}$.
    
    We will use the following $\epsilon$-dependent form of the Park--Pham theorem, due to Bell \cite{bell-park-pham}.  We write logarithms in base $2$.
    
    \begin{theorem}[Bell's $\epsilon$-dependent Park--Pham theorem \cite{bell-park-pham}]
    \label{bell optimal bound theorem}
        Let $X$ be a finite set, let $\mathcal{F}\subseteq 2^X$ be a nontrivial upper set, and let $\epsilon\in(0,1)$.  If $p>48 q(\mathcal{F})\log_2\left(\frac{\ell(\mathcal{F})}{\epsilon}\right)$, then $\mathbb P(X_p\in\mathcal{F})>1-\epsilon$.
    \end{theorem}

    This formulation follows from Bell's original statement \cite{bell-park-pham} by applying it with any $q'>q(\mathcal F)$ sufficiently close to $q(\mathcal F)$.
    
    The purpose of this note is to record a conditional consequence of Theorem~\ref{bell optimal bound theorem}.  Suppose that one does not sample from the full Boolean lattice $2^X$, but instead samples from $\mu_p$ conditioned on some event $B\subseteq 2^X$.  We ask when the Park--Pham lower bound for an upper set $\mathcal{F}$ transfers to the conditional probability $\mathbb P(X_p\in\mathcal{F}\mid X_p\in B)$.  The transfer mechanism is elementary, but useful: the precise hypothesis needed for conditioning to preserve the Park--Pham lower bound is positive correlation with the target upper set.  This observation gives a convenient way to move Park--Pham estimates from the full Boolean lattice to conditional models, including monotone conditioning events, independent-coordinate restrictions, and finite-poset models embedded in Boolean lattices.
    
    This framing is deliberately conditional.  Not every conditioning event $B$ is favorable: if $B$ is negatively correlated with $\mathcal{F}$, or if $B$ is disjoint from $\mathcal{F}$, then no such transfer can be expected.  However, the positive-correlation condition is checkable in several natural cases.  If $B$ is itself an upper set, then Harris's inequality gives the required positive correlation.  If $B$ depends on coordinates disjoint from those determining $\mathcal{F}$, then independence gives the same conclusion, even when $B$ is not monotone.
    
    We also record a poset-valued formulation.  If $P$ is a nonempty finite poset, $U\subseteq P$ is an upper set, and $f:P\to 2^X$ is an order embedding, then conditioning $X_p$ on the image $f(P)$ induces a $P$-valued random variable.  Under the corresponding positive-correlation hypothesis, the Park--Pham bound for the Boolean upper set generated by $f(U)$ yields a conditional threshold estimate for $U$.

    Thus, the note is organizational as much as technical: we isolate the exact positive-correlation condition under which the $\epsilon$-dependent Park--Pham estimate survives conditioning.  This gives a reusable criterion for conditional threshold models, especially in settings where the natural probability space is not the full Boolean lattice but a conditioned subspace or an embedded finite poset.
    
\section{Conditional transfer under positive correlation}

    We begin by isolating the elementary transfer mechanism.
    
    \begin{definition}
        Let $\mathcal A,\mathcal B\subseteq 2^X$.  We say that $\mathcal A$ and $\mathcal B$ are \textit{positively correlated at $p$} if
        \[
            \mu_p(\mathcal A\cap\mathcal B) \geq \mu_p(\mathcal A)\mu_p(\mathcal B).
        \]
    \end{definition}
    
    \begin{theorem}[Conditional Park--Pham transfer]
    \label{positive correlation theorem}
        Let $X$ be a finite set, let $p\in [0,1]$, let $\mathcal{F}\subseteq 2^X$ be a nontrivial upper set, and let $B\subseteq 2^X$ satisfy $\mu_p(B)>0$.  Suppose that $\mathcal{F}$ and $B$ are positively correlated at $p$.  If $\epsilon\in(0,1)$ and $p>48 q(\mathcal{F})\log_2\left(\frac{\ell(\mathcal{F})}{\epsilon}\right)$, then $\mathbb P(X_p\in\mathcal{F}\mid X_p\in B)>1-\epsilon$.
    \end{theorem}
    
    \begin{proof}
        By positive correlation,
        \[
            \mathbb P(X_p\in\mathcal{F}\mid X_p\in B) = \frac{\mu_p(\mathcal{F}\cap B)}{\mu_p(B)} \geq \mu_p(\mathcal{F}).
        \]
        Theorem~\ref{bell optimal bound theorem} gives $\mu_p(\mathcal{F})>1-\epsilon$ under the displayed hypothesis on $p$.
    \end{proof}
    
    \begin{remark}
        The positive-correlation hypothesis in Theorem~\ref{positive correlation theorem}
        is pointwise in $p$.  Thus, in applications, the conclusion holds only for
        values of $p$ for which both the Park--Pham lower bound and the correlation
        inequality
        \[
            \mu_p(\mathcal{F}\cap B)\geq \mu_p(\mathcal{F})\mu_p(B)
        \]
        hold simultaneously.  The sufficient conditions given below, namely Harris's
        inequality for increasing conditioning events and independence for disjoint
        coordinate restrictions, have the advantage of holding uniformly for all
        $p\in(0,1)$.
    \end{remark}

\section{Sources of positive correlation}

    The simplest source of positive correlation is monotonicity.  When both the target event and the conditioning event are upper sets, Harris's inequality shows that the conditioning can only increase the probability of the target event.  The more flexible point of Theorem~\ref{positive correlation theorem} is that the same transfer also applies to nonmonotone conditioning events whenever positive correlation can be verified by some other means.
    
    \begin{theorem}[Harris's inequality \cite{harris}]
    \label{harris inequality}
        If $\mathcal A,\mathcal B\subseteq 2^X$ are upper sets, then $\mu_p(\mathcal A\cap\mathcal B) \geq \mu_p(\mathcal A)\mu_p(\mathcal B)$ for every $p\in[0,1]$.
    \end{theorem}
    
    \begin{remark}
        When both $\mathcal F$ and $B$ are upper sets, Harris's inequality gives
        \[
            \mathbb P(X_p\in\mathcal F\mid X_p\in B)\geq \mu_p(\mathcal F),
        \]
        so Bell's theorem (i.e. Theorem~\ref{bell optimal bound theorem}) immediately yields the same conditional bound.  The point of Theorem~\ref{positive correlation theorem} is that the same argument also applies when positive correlation is obtained by other means, including nonmonotone conditioning events.
    \end{remark}
    
    The positive-correlation hypothesis is not limited to monotone conditioning events.  The following simple case shows that the same conclusion can hold for nonmonotone restrictions.
    
    \begin{proposition}[Independent-coordinate conditioning]
    \label{independent coordinate proposition}
        Let $X=Y\sqcup Z$ be a disjoint union.  Let $\mathcal G\subseteq 2^Y$ be a nontrivial upper set, and define $\mathcal{F}=\{S\subseteq X:S\cap Y\in\mathcal G\}$.  Let $\mathcal H\subseteq 2^Z$ be nonempty, and define $B=\{S\subseteq X:S\cap Z\in\mathcal H\}$.  Then, $\mathcal{F}$ and $B$ are independent under $\mu_p$ for every $p\in(0,1)$.  In particular, if $\epsilon\in(0,1)$ and $p>48 q(\mathcal{F})\log_2\left(\frac{\ell(\mathcal{F})}{\epsilon}\right)$, then $\mathbb P(X_p\in\mathcal{F}\mid X_p\in B)>1-\epsilon$.
    \end{proposition}
    
    \begin{proof}
        The event $\mathcal{F}$ depends only on the coordinates in $Y$, while $B$ depends only on the coordinates in $Z$.  Since $\mu_p$ is a product measure and $Y\cap Z=\emptyset$, the two events are independent and $\mu_p(\mathcal{F}\cap B)=\mu_p(\mathcal{F})\mu_p(B)$.  The desired conditional bound follows from Theorem~\ref{positive correlation theorem}.
        \end{proof}
        
        \begin{remark}
        The preceding proposition includes many nonmonotone conditioning events.  For instance, one may take
        \[
            \mathcal H=\{T\subseteq Z: |T|\text{ is even}\},
        \]
        so that $B$ is the event that an even number of $Z$-coordinates are selected.  This event is not an upper set, but it is still independent of every event determined only by the $Y$-coordinates.
    \end{remark}

\section{Poset-valued formulations and Boolean lifts}

    We now translate the preceding conditional bound into a finite-poset setting. First, we record a general, flexible formulation where we embed a poset into a boolean lattice such that the image of the poset is positively correlated with the upper set $\mathcal{F}$.

    \begin{theorem}[Boolean lifts of poset upper sets]\label{boolean}
        Let $p\in(0,1)$, let $P$ be a nonempty finite poset, let $f:P\to 2^X$ be an injection, and let $U\subseteq P$.  Suppose that $\mathcal{F}\subseteq 2^X$ is a nontrivial upper set satisfying $\mathcal{F}\cap f(P)=f(U)$.  Let $Y_p$ be the $P$-valued random variable obtained by conditioning $X_p$ on $f(P)$ and pulling back along $f$. If $\epsilon\in(0,1)$, $\mathcal{F}$ and $f(P)$ are positively correlated at $p$, and $p>48q(\mathcal{F})\log_2\left(\frac{\ell(\mathcal{F})}{\epsilon}\right)$, then $\mathbb P(Y_p\in U)>1-\epsilon$.
    \end{theorem}
    
    \begin{proof}
        Since $\mathcal{F}\cap f(P)=f(U)$, we have
        \[
            \mathbb P(Y_p\in U) = \frac{\mu_p(f(U))}{\mu_p(f(P))} = \frac{\mu_p(\mathcal{F}\cap f(P))}{\mu_p(f(P))}.
        \]
        By positive correlation, this is at least $\mu_p(\mathcal{F})$.  Bell's
        $\epsilon$-dependent Park--Pham theorem gives the result.
    \end{proof} 
    A particularly useful way to produce Boolean lifts is through order embeddings.  Let $P$ be a nonempty finite poset and recall that an injective map $f:P\to 2^X$ is an \textit{order embedding} if, for all $x,y\in P$, $x\leq_P y \quad\Longleftrightarrow\quad f(x)\subseteq f(y)$.    If $f:P\to 2^X$ is an order embedding and $U\subseteq P$ is an upper set, then the Boolean upper set $\langle f(U)\rangle$ cuts out $f(U)$ inside $f(P)$.
     
    \begin{lemma}
    \label{poset intersection lemma}
        Let $P$ be a nonempty finite poset, let $U\subseteq P$ be an upper set, and let $f:P\to 2^X$ be an order embedding.  Put $\mathcal{F}=\langle f(U)\rangle$.  Then, $\mathcal{F}\cap f(P)=f(U)$.
    \end{lemma}
    
    \begin{proof}
        The inclusion $f(U)\subseteq \mathcal{F}\cap f(P)$ is immediate.  Conversely, suppose $T\in \mathcal{F}\cap f(P)$.  Then $T=f(y)$ for some $y\in P$, and since $T\in\langle f(U)\rangle$, there exists $x\in U$ such that $f(x)\subseteq f(y)$.  Since $f$ is an order embedding, $x\leq_P y$.  Because $U$ is an upper set, $y\in U$.  Hence $T=f(y)\in f(U)$.
    \end{proof}
    If $f:P\to 2^X$ is an injection and $p\in(0,1)$, then conditioning $\mu_p$ on the event $X_p\in f(P)$ induces a probability measure $\nu_p$ on $P$ by
    \[
        \nu_p(A)
        :=
        \frac{\mu_p(f(A))}{\mu_p(f(P))}
        \qquad (A\subseteq P).
    \]
    Equivalently, if $Y_p$ has distribution $\nu_p$, then $\mathbb{P}(Y_p\in A) = \mathbb{P}(X_p\in f(A)\mid X_p\in f(P))$.
    
    \begin{corollary}[Embedded finite posets under positive correlation]
    \label{poset positive correlation theorem}
         Let $p\in(0,1)$, let $P$ be a nonempty finite poset, let $U\subseteq P$ be a nonempty upper set, and let $f:P\to 2^X$ be an order embedding.  Put $\mathcal{F}=\langle f(U)\rangle$, set $B=f(P)$, and assume that $\mathcal{F}$ is nontrivial.  Suppose that $\mathcal{F}$ and $B$ are positively correlated at $p$.  If $\epsilon\in(0,1)$ and $p>48 q(\mathcal{F})\log_2\left(\frac{\ell(\mathcal{F})}{\epsilon}\right)$, then $\mathbb P(Y_p\in U)>1-\epsilon$.
    \end{corollary}
    
    \begin{proof}
        By Lemma~\ref{poset intersection lemma}, $f(U)=\mathcal{F}\cap f(P)=\mathcal{F}\cap B$.  Therefore,
        \[
            \mathbb P(Y_p\in U)
            =
            \frac{\mu_p(f(U))}{\mu_p(f(P))}
            =
            \frac{\mu_p(\mathcal{F}\cap B)}{\mu_p(B)}.
        \]
        Since $\mathcal{F}$ and $B$ are positively correlated at $p$, this is at least $\mu_p(\mathcal{F})$.  Theorem~\ref{bell optimal bound theorem} now gives $\mu_p(\mathcal{F})>1-\epsilon$ under the displayed hypothesis.
    \end{proof}
    
    The following corollary gives a structural sufficient condition for the positive-correlation hypothesis in Corollary~\ref{poset positive correlation theorem}.
    
    \begin{corollary}[Upper-image embedded posets]
    \label{upper image poset corollary}
        Let $P$ be a nonempty finite poset, let $U\subseteq P$ be a nonempty upper set, and let $f:P\to 2^X$ be an order embedding such that $f(P)$ is an upper set in $2^X$.  Put $\mathcal{F}=\langle f(U)\rangle$, and assume that $\mathcal{F}$ is nontrivial.  If $p\in (0,1)$, $\epsilon\in(0,1)$, and $p>48 q(\mathcal{F})\log_2\left(\frac{\ell(\mathcal{F})}{\epsilon}\right)$, then $\mathbb P(Y_p\in U)>1-\epsilon$.
    \end{corollary}
    
    \begin{proof}
        The family $\mathcal{F}=\langle f(U)\rangle$ is an upper set by construction, and $f(P)$ is an upper set by assumption.  Harris's inequality implies that $\mathcal{F}$ and $f(P)$ are positively correlated for every $p\in[0,1]$.  The result follows from Corollary~\ref{poset positive correlation theorem}.
    \end{proof}
    
    \begin{remark}
        The assumption that $f(P)$ is an upper set is restrictive, but it is a clean, checkable, condition.  More generally, Corollary~\ref{poset positive correlation theorem} applies whenever the image $f(P)$ is positively correlated with the Boolean upper set generated by $f(U)$.  In applications where $f(P)$ is not monotone, this correlation condition may still hold by independence, symmetry, or direct computation.
    \end{remark}

\section{Examples}
    
    We give examples illustrating the preceding mechanisms and the permutation-profile formulation.    
    
    \begin{example}[A nonmonotone independent conditioning event]
        Let $X=Y\sqcup Z$, and let $\mathcal G\subseteq 2^Y$ be a nontrivial upper set.  Define $\mathcal{F} = \{S\subseteq X:S\cap Y\in\mathcal G\}$ and condition on the event $B=\{S\subseteq X: |S\cap Z|\text{ is even}\}$.  The event $B$ is generally not an upper set.  However, $\mathcal{F}$ depends only on the $Y$-coordinates and $B$ depends only on the $Z$-coordinates, so the two events are independent under $\mu_p$.  Hence, Proposition~\ref{independent coordinate proposition} gives $\mathbb P(X_p\in\mathcal{F}\mid X_p\in B)>1-\epsilon$ whenever $p>48 q(\mathcal{F})\log_2\left(\frac{\ell(\mathcal{F})}{\epsilon}\right)$. Thus the Park--Pham bound transfers unchanged even though the conditioning event is not monotone.
    \end{example}
    
    \begin{example}[Permutation pattern classes]\label{permutation-example}
        Let $C$ be a permutation class, and let $C_{\leq n}$ be the nonempty finite poset of permutations in $C$ of length at most $n$, ordered by pattern containment.  Define
        \[
            \Phi(\sigma)=\{\tau\in C_{\leq n}:\tau\preceq\sigma\}
        \]
        for $\sigma\in C_{\leq n}$.  Then $\Phi:C_{\leq n}\to 2^{C_{\leq n}}$ is an order embedding.  Let $\mathcal U\subseteq C_{\leq n}$ be a nontrivial upper set, and let $M$ be the set of minimal elements of $\mathcal U$.  Define
        \[
            \mathcal E_M=\{T\subseteq C_{\leq n}:T\cap M\neq\emptyset\}.
        \]
        Then $\mathcal E_M$ is an upper set in the Boolean lattice $2^{C_{\leq n}}$, and $\mathcal E_M\cap \Phi(C_{\leq n})=\Phi(\mathcal U)$.  Indeed, for $\sigma\in C_{\leq n}$,
        \[
            \Phi(\sigma)\in\mathcal E_M
            \quad\Longleftrightarrow\quad
            \sigma\text{ contains some }\pi\in M
            \quad\Longleftrightarrow\quad
            \sigma\in\mathcal U.
        \]

        Write  $m=|M|$.  The minimal elements of $\mathcal E_M$ are precisely the singletons $\{\pi\}$ with $\pi\in M$, so $\ell(\mathcal E_M)=1$.  Moreover, $\mathcal E_M$ is covered by these $m$ singletons, giving $q(\mathcal E_M)\geq 1/(2m)$.  Conversely, if $\mathcal G$ covers $\mathcal E_M$, then for each $\pi\in M$, the singleton $\{\pi\}\in\mathcal E_M$ must contain some $G\in\mathcal G$. Since the only subsets of $\{\pi\}$ are $\emptyset$ and $\{\pi\}$, and including $\emptyset$ would contribute $1$ to the covering sum, any $p$-small cover must contain $\{\pi\}$ for each $\pi\in M$. Hence the covering sum is at least $mp$, so $q(\mathcal E_M)=1/(2m)$.
        
        Therefore, if $\mathcal E_M$ and $\Phi(C_{\leq n})$ are positively correlated at $p$, then $\mathbb P(Y_p\in\mathcal U)>1-\epsilon$ whenever 
        \[
            p>48q(\mathcal{E}_M)\log_2\left(\frac{\ell(\mathcal{E}_M)}{\epsilon}\right)=\frac{24}{m}\log_2\left(\frac{1}{\epsilon}\right).
        \]
    \end{example}
    
    \begin{remark}
        In the permutation-pattern setting, the positive-correlation hypothesis is a genuine additional condition.  Determining when this holds for natural permutation classes appears to be an interesting problem in its own right.  The image $\Phi(C_{\leq n})$ is typically far from an upper set in $2^{C_{\leq n}}$, so Harris's inequality does not apply directly.  Thus the conditional Park--Pham transfer reduces the problem to understanding when the set of valid pattern profiles is positively correlated with the Boolean event of hitting the minimal pattern antichain $M$, rather than the stronger condition of containing an entire prescribed pattern profile.
    \end{remark}

    The preceding example leaves the positive-correlation condition as an additional hypothesis.  The next example shows that this condition is nonempty by verifying it in a natural permutation class.

    \begin{example}[Stack-sortable permutations]
        Let $C=\operatorname{Av}(231)$ be the class of one-stack-sortable permutations, where $\operatorname{Av}(231)$ denotes the permutations avoiding the pattern $231$.  Thus $|C_r|=\operatorname{Cat}_r$, the $r$-th Catalan number \cite{knuth}.  We take the convention that $C_0$ consists of the empty permutation. Fix $1\leq k\leq n$, and let $C_{\le n}= C_0\cup\hdots\cup C_{n}$ be ordered by pattern containment.  Define
        \[
            \mathcal U_{k,n}=\{\sigma\in C_{\le n}:|\sigma|\ge k\}.
        \]
        Then $\mathcal U_{k,n}$ is an upper set, and its minimal elements are exactly $M=C_k$.  In particular, $m=|M|=\operatorname{Cat}_k$.  Let $\Phi(\sigma)$ and $\mathcal{E}_M$ be as before, so that $\mathcal{E}_M\cap \Phi(C_{\le n})=\Phi(\mathcal{U}_{k,n})$.
        
        At $p=1/2$, the measure on $2^{C_{\le n}}$ is uniform.  Hence, we get 
        \[
            \mathbb P(Y_{1/2}\in\mathcal U_{k,n}) = \frac{\mu_{1/2}\left(\mathcal{E}_M \cap \Phi(C_{\leq n})\right)}{\mu_{1/2}\big(\Phi(C_{\leq n})\big)}= \frac{|\mathcal U_{k,n}|}{|C_{\le n}|} = \frac{\sum_{r=k}^n \operatorname{Cat}_r}{\sum_{r=0}^n \operatorname{Cat}_r}
        \]
        and $\mu_{1/2}(\mathcal E_M)=1-\left(1-\frac{1}{2}\right)^{\operatorname{Cat}_k}=1-2^{-\operatorname{Cat}_k}$.  Therefore, $\mathcal E_M$ and $\Phi(C_{\le n})$ are positively correlated at $p=1/2$ whenever
        \begin{align*}
            \frac{\mu_{1/2}\left(\mathcal{E}_M \cap \Phi(C_{\leq n})\right)}{\mu_{1/2}\big(\Phi(C_{\leq n})\big)} & \geq \mu_{1/2}(\mathcal E_M)\\
            \frac{\sum_{r=k}^n \operatorname{Cat}_r}{\sum_{r=0}^n \operatorname{Cat}_r} &\geq 1-2^{-\operatorname{Cat}_k} \\
            2^{-\operatorname{Cat}_k}\sum_{r=0}^n \operatorname{Cat}_r  &\geq \left(\sum_{r=0}^n \operatorname{Cat}_r\right)-\left(\sum_{r=k}^n \operatorname{Cat}_r\right) = \sum_{r=0}^{k-1} \operatorname{Cat}_r \\
        \end{align*}
        For each fixed $k$, this holds for all sufficiently large $n$.  Since $q(\mathcal E_M)=(2\operatorname{Cat}_k)^{-1}$ and $\ell(\mathcal E_M)=1$, the conditional Park--Pham transfer result (Theorem~\ref{boolean}) gives $\mathbb P(Y_{1/2}\in\mathcal U_{k,n})>1-\epsilon$ whenever $\epsilon>2^{-\operatorname{Cat}_k/48}$ and $n$ is sufficiently large.
    \end{example}
    
    \begin{remark}
        The preceding example is meant as a calibration example rather than a deep application: the conclusion can also be obtained directly from the Catalan enumeration.  Its purpose is to show that the positive-correlation hypothesis in the permutation-profile formulation is not empty, even for a natural non-Boolean permutation poset.  More interesting applications would require verifying the same correlation inequality for upper sets defined by a proper subcollection of patterns, rather than by all patterns of a fixed length.
    \end{remark}
\section{Concluding remarks}

    The preceding results show that the $\epsilon$-dependent Park--Pham theorem admits a simple conditional form under positive correlation.  The point is not that every conditioning event preserves a Park--Pham-type bound; rather, the positive-correlation condition is exactly the hypothesis ensuring that conditioning does not reduce the probability of the target upper set.  Harris's inequality supplies one general source of such conditioning events, while independent-coordinate restrictions supply another, including nonmonotone examples.
    
    In light of Example~\ref{permutation-example}, it would also be interesting to determine for which permutation classes $C$, integers $n$, and upper sets $\mathcal U\subseteq C_{\leq n}$ with minimal antichain $M$, the inequality
    \[
        \mu_p(\mathcal E_M\cap \Phi(C_{\leq n})) \geq
        \mu_p(\mathcal E_M)\mu_p(\Phi(C_{\leq n}))
    \]
    holds.  Even special cases would give conditional Park--Pham bounds for pattern-containment properties, or equivalently for complements of pattern-avoidance properties.

\end{document}